\documentclass[12pt,reqno]{amsart}
\usepackage{amsmath,amssymb,amsfonts,amscd,latexsym,amsthm,mathrsfs,mathtools,verbatim,cite,comment}
\usepackage[unicode]{hyperref}
\let\ol\overline

\begin{document}

\title[First memoir on the asymptotics of certain infinite products]{First memoir on the asymptotics\\of certain infinite products}

\author{Wadim Zudilin}
\address{Department of Mathematics, IMAPP, Radboud University, PO Box 9010, 6500~GL Nijmegen, The Netherlands}
\urladdr{https://www.math.ru.nl/~wzudilin/}

\dedicatory{To Mourad Ismail, with admiration and warm wishes}

\date{25 November 2024}

\subjclass[2020]{Primary 11P84; Secondary 05A15, 11F03, 11N37}
\keywords{Rogers--Ramanujan identities; modular function; modular equation; asymptotics; Kanade--Russell conjectures}

\begin{abstract}
The product sides of the Rogers--Ramanujan identities and alike often appear to be `transparently modular' (functions). The old work by Rogers (1894) and recent work by Rosengren make use (somewhat implicitly) of this fact for proving the identities with the help of underlying modular equations\,---\,the main challenge is verifying the latter for the sum sides.
Here we speculate on the potentials of using the asymptotics of such $q$-identities or their finite versions for proving them.
\end{abstract}

\thanks{The work was supported by the Max Planck Institute for Mathematics (Bonn).}

\maketitle

This note is inspired by several developments around the famous Rogers--Ramanujan identities including the old memoir \cite{Ro94} of L.\,J.~Rogers, which made an influence on the present title.
It would be also not mistaken to acknowledge that my personal interest in, taste for and education on $q$-series and Rogers--Ramanujan identities (in particular) were tremendously influenced by the work of Mourad Ismail.
Particular examples of his fine work on the latter subject, again through my personal views, are \cite{GIS99,IS06,IZ18}.
I find it appropriate to dedicate this piece to Mourad on the occasion of his round birthday.

To prepare the stage, we list the standard $q$-notation (when $q$ is viewed as a complex parameter we assume it to satisfy $|q|<1$): the $q$-Pochhammer symbol is given by
\[
(a;q)_k=\prod_{j=0}^{k-1}(1-aq^j)
\]
and its multiple version is $(a_1,\dots,a_s;q)_k=(a_1;q)_k\dotsb(a_s;q)_k$.
These are used below for $k$ a non-negative integer as well as for $k=\infty$.

\section{Rogers--Ramanujan and Rogers--Selberg identities}
\label{sec:1}

One proof of the celebrated Rogers--Ramanujan identities 
\begin{align*}
G(q)=\sum_{k=0}^\infty\frac{q^{k^2}}{(q;q)_k}
=\frac 1{(q,q^4;q^5)_\infty},
\quad
H(q)=\sum_{k=0}^\infty\frac{q^{k(k+1)}}{(q;q)_k}
=\frac 1{(q^2,q^3;q^5)_\infty},
\end{align*}
based on the original methodology of Rogers \cite{Ro94} was recently given by Rosengren~\cite{Ro24}.
It is based on the `self-replicating' equations satisfied by (both sum and product sides of) the two functions:
\begin{equation} \label{RR}
\begin{aligned}
G(q)&=\frac{(q^8;q^8)_\infty}{(q^2;q^2)_\infty}\big(qH(-q^4)+G(q^{16})\big),
\\
H(q)&=\frac{(q^8;q^8)_\infty}{(q^2;q^2)_\infty}\big(G(-q^4)+q^3H(q^{16})\big).
\end{aligned}
\end{equation}
Similar machinery was used by Rosengren in \cite{Ro25} to prove the identities
\begin{align*}
A(q)&=(-q;q)_\infty\sum_{k=0}^\infty\frac{q^{2k^2}}{(q^4;q^4)_k(-q;q^2)_k}
=\frac{1}{(q,q^2,q^5,q^6;q^7)_\infty},
\\
B(q)&=(-q;q)_\infty\sum_{k=0}^\infty\frac{q^{2k(k+1)}}{(q^4;q^4)_k(-q;q^2)_k}
=\frac{1}{(q,q^3,q^4,q^6;q^7)_\infty},
\\
C(q)&=(-q;q)_\infty\sum_{k=0}^\infty\frac{q^{2k(k+1)}}{(q^4;q^4)_k(-q;q^2)_{k+1}}
=\frac{1}{(q^2,q^3,q^4,q^5;q^7)_\infty}
\end{align*}
by showing that
\begin{equation} \label{septic}
\begin{gathered}
A(q)=\frac 1{(q^2;q^4)_\infty^2}\big(qC(-q^2)+A(q^8)\big),
\quad
B(q)=\frac 1{(q^2;q^4)_\infty^2}\big(A(-q^2)+qB(q^8)\big),
\\
C(q)=\frac 1{(q^2;q^4)_\infty^2}\big(B(-q^2)+q^3C(q^8)\big),
\end{gathered}
\end{equation}
for both the sum and product sides.

The first thing to point out is that the product sides of all the identities above are \emph{modular} functions after an appropriate normalisation (which makes identities \eqref{RR} and \eqref{septic} particular examples of modular equations for participating modular functions).
Namely, for $b>a>0$ integers, the product $(q^a,q^{b-a};q^b)_\infty$ transforms into a modular function
\[
q^{bB_2(a/b)/2}(q^a,q^{b-a};q^b)_\infty,
\]
where $B_2(t)=t^2-t+1/6$ is the second Bernoulli polynomial;
also $q^{b/24}(q^b;q^b)_\infty$ is a (weight $1/2$) modular function (a `scaled' Dedekind eta function).
For example, if we write
\begin{gather*}
\hat A(q)=q^{-1/42}A(q)=\frac{q^{-1/42}}{(q,q^2,q^5,q^6;q^7)_\infty},
\quad
\hat B(q)=q^{5/42}B(q)=\frac{q^{5/42}}{(q,q^3,q^4,q^6;q^7)_\infty},
\\
\hat C(q)=q^{17/42}C(q)=\frac{q^{17/42}}{(q^2,q^3,q^4,q^5;q^7)_\infty},
\end{gather*}
then equations \eqref{septic} assume a symmetric form
\begin{gather*}
\hat A(q)=\frac{q^{1/6}}{(q^2;q^4)_\infty^2}\big(\pm\hat C(-q^2)+\hat A(q^8)\big),
\quad
\hat B(q)=\frac{q^{1/6}}{(q^2;q^4)_\infty^2}\big(\pm\hat A(-q^2)+\hat B(q^8)\big),
\\
\hat C(q)=\frac{q^{1/6}}{(q^2;q^4)_\infty^2}\big(\pm\hat B(-q^2)+\hat C(q^8)\big)
\end{gather*}
(we ignore the choice of sign), in which $q^{1/6}$ corresponds to the modular normalisation of the eta-type product $1/(q^2;q^4)_\infty^2$.
In fact, the sign variation leads to companion identities
\begin{align*}
A(q)&=\frac{(-q,-q^5,-q^9,-q^{13};q^{14})_\infty\cdot\big(-qC(-q^2)+A(q^8)\big)}{(q^2;q^4)_\infty^2(q,q^5,q^9,q^{13};q^{14})_\infty},
\\
B(q)&=\frac{(-q,-q^3,-q^{11},-q^{13};q^{14})_\infty\cdot\big(A(-q^2)-qB(q^8)\big)}{(q^2;q^4)_\infty^2(q,q^3,q^{11},q^{13};q^{14})_\infty},
\\
C(q)&=\frac{(-q^3,-q^5,-q^9,-q^{11};q^{14})_\infty\cdot\big(B(-q^2)-q^3C(q^8)\big)}{(q^2;q^4)_\infty^2(q^3,q^5,q^9,q^{11};q^{14})_\infty}
\end{align*}
of \eqref{septic}, and we also have similar closed forms for
\begin{align*}
G(q)
&=\frac{-qH(-q^4)+G(q^{16})}{(q^4;q^8)_\infty(q,q^9;q^{10})_\infty^2(q^6,q^{10},q^{14};q^{20})_\infty},
\\
H(q)
&=\frac{G(-q^4)-q^3H(q^{16})}{(q^4;q^8)_\infty(q^3,q^7;q^{10})_\infty^2(q^2,q^{10},q^{18};q^{20})_\infty}
\end{align*}
complementing \eqref{RR}.

Notice however that these alternative modular equations for $A(q),B(q),C(q)$ and for $G(q),H(q)$ seem to be harder to establish directly for the sum sides of the corresponding Rogers--Ramanujan-type identities.

\section{Kanade--Russell modulo 9 identities}
\label{sec:2}

One can also perform similar modular normalisation of the product sides of the Kanade--Russell modulo 9 conjectural identities \cite{KR15,UZ21}
\begin{align*}
\sum_{m,n\ge0} \frac{q^{m^2+3mn+3n^2}}{(q;q)_m(q^3;q^3)_n}
\overset?=K_1(q)
&=\frac{1}{(q,q^3,q^6,q^8;q^9)_\infty},
\\
\sum_{m,n\ge0} \frac{q^{m^2+3mn+3n^2+m+3n}}{(q;q)_m(q^3;q^3)_n}
\overset?=K_2(q)
&=\frac{1}{(q^2,q^3,q^6,q^7;q^9)_\infty},
\\
\sum_{m,n\ge0} \frac{q^{m^2+3mn+3n^2+2m+3n}}{(q;q)_m(q^3;q^3)_n}
\overset?=K_3(q)
&=\frac{1}{(q^3,q^4,q^5,q^6;q^9)_\infty},
\end{align*}
to check (routinely!) that these product sides satisfy the equations
\begin{align*}
K_1(q)
&=\frac{qK_2(-q^2)+K_1(q^8)}{(q^2;q^4)_\infty(q^6;q^{12})_\infty(q^3,q^{15};q^{18})_\infty(-q^5,-q^{13};q^{18})_\infty},
\\
K_2(q)
&=\frac{K_3(-q^2)+qK_2(q^8)}{(q^2;q^4)_\infty(q^6;q^{12})_\infty(q^3,q^{15};q^{18})_\infty(-q,-q^{17};q^{18})_\infty},
\\
K_3(q)
&=\frac{K_1(-q^2)+q^5K_3(q^8)}{(q^2;q^4)_\infty(q^6;q^{12})_\infty(q^3,q^{15};q^{18})_\infty(-q^7,-q^{11};q^{18})_\infty}.
\end{align*}
They also satisfy
\begin{align*}
K_1(q)
&=\frac{(-q,-q^{17};q^{18})_\infty\cdot\big(-qK_2(-q^2)+K_1(q^8)\big)}{(q^2;q^4)_\infty(q^6;q^{12})_\infty(q^3,q^{15};q^{18})_\infty(q,q^5,q^{13},q^{17};q^{18})_\infty},
\\
K_2(q)
&=\frac{(-q^7,-q^{11};q^{18})_\infty\cdot\big(K_3(-q^2)-qK_2(q^8)\big)}{(q^2;q^4)_\infty(q^6;q^{12})_\infty(q^3,q^{15};q^{18})_\infty(q,q^7,q^{11},q^{17};q^{18})_\infty},
\\
K_3(q)
&=\frac{(-q^5,-q^{13};q^{18})_\infty\cdot\big(K_1(-q^2)-q^5K_3(q^8)\big)}{(q^2;q^4)_\infty(q^6;q^{12})_\infty(q^3,q^{15};q^{18})_\infty(q^5,q^7,q^{11},q^{13};q^{18})_\infty}.
\end{align*}
The techniques in \cite{Ro24,Ro25} are hardly applicable to the sum sides of either identities for $K_1(q),\allowbreak K_2(q),K_3(q)$: this time we deal with double summations.

\section{Andrews--Gordon identities}
\label{sec:3}

The Rogers--Ramanujan identities are first entries in a general family of the Andrews--Gordon identities, which feature for each $k\ge2$ the product sides.
Similarly, the product parts of modulo 7 Andrews--Gordon identities,
\[
P_i(q)=\frac{(q^i,q^{2k+1-i},q^{2k+1};q^{2k+1})_\infty}{(q;q)_\infty}, \quad\text{where} \; i=1,\dots,k,
\]
while the sum sides correspond to $(k-1)$-fold summations.
For $k=3$, the eta-type products coincide with those for $A(q),B(q),C(q)$ which we have already treated in Section~\ref{sec:1}.
Similar treatment for the sum sides is not known.

Unfortunately, there seems to be no simple-looking extension of the story to modulo~9,~13 (and higher).
For modulo~11 and the products
\[
P_i(q)=\frac{(q^i,q^{11-i},q^{11};q^{11})_\infty}{(q;q)_\infty}, \quad\text{where} \; i=1,\dots,5,
\]
we find some `surrogate' versions:
\begin{align*}
P_1(q^2)&=\frac{P_5(-q^2)-q^6P_3(q^8)}{(q^2;q^4)_\infty^2(q^{18},q^{26};q^{44})_\infty},
\\
P_2(q^2)&=\frac{-q^2P_1(-q^2)+P_5(q^8)}{(q^2;q^4)_\infty^2(q^{14},q^{30};q^{44})_\infty},
\\
P_3(q^2)&=\frac{P_4(-q^2)-q^4P_2(q^8)}{(q^2;q^4)_\infty^2(q^{10},q^{34};q^{44})_\infty},
\\
P_4(q^2)&=\frac{P_2(-q^2)-q^6P_1(q^8)}{(q^2;q^4)_\infty^2(q^6,q^{38};q^{44})_\infty},
\\
P_5(q^2)&=q^2\,\frac{-P_3(-q^2)+P_4(q^8)}{(q^2;q^4)_\infty^2(q^2,q^{42};q^{44})_\infty};
\end{align*}
and there are no alternative plus-sign companions.

\section{Asymmetric Kanade--Russell identities}
\label{sec:4}

There are more Kanade--Russell modulo 9 identities \cite{KR15,UZ21} which are not modular functions.
They include
\begin{align*}
\sum_{m,n\ge0} \frac{q^{m^2+3mn+3n^2+m+2n}}{(q;q)_m (q^3;q^3)_n}
\overset?=K_4(q)
&=\frac{1}{(q^2,q^3,q^5,q^8;q^9)_\infty}
=\frac{(q^9;q^9)_\infty(\omega q^2,\ol\omega q^2;q^3)_\infty}{(q^3;q^3)_\infty},
\\
\sum_{m,n\ge0} \frac{q^{m^2+3mn+3n^2-m+n}(1-q^{2m})}{(q;q)_m (q^3;q^3)_n}
\overset?=K_5(q)
&=\frac{1}{(q,q^4,q^6,q^7;q^9)_\infty}
=\frac{(q^9;q^9)_\infty(\omega q,\ol\omega q;q^3)_\infty}{(q^3;q^3)_\infty},
\end{align*}
where the sum side of the last identity was simplified from its original by Hickerson \cite{Hi21}, who also complemented the set with the identities
\begin{equation} \label{prod}
\begin{aligned}
\sum_{m,n\ge0} \frac{q^{m^2+3mn+3n^2+m+n}(1-\omega q^{m+3n+1})}{(q;q)_m(q^3;q^3)_n}
\overset?=K_6(q)
&=\frac{(q^6;q^9)_\infty(\omega q,\ol\omega q^3;q^3)_\infty}{(q^2;q^3)_\infty},
\\
\sum_{m,n\ge0} \frac{q^{m^2+3mn+3n^2+2n}(1-\omega q^{3m+3n+2})}{(q;q)_m(q^3;q^3)_n}
\overset?=K_7(q)
&=\frac{(q^3;q^9)_\infty(\omega q^2,\ol\omega q^3;q^3)_\infty}{(q;q^3)_\infty},
\end{aligned}
\end{equation}
with $\omega$ a primitive cubic root of unity and $\ol\omega$ its conjugate,
and their conjugations (so that there are four additional conjectural identities in total).
It is hard to expect any modular-type functional equations for these new instances, though the product sides may possess some modular-like behaviour.

Consider a real number $a$ from the interval $0<a\le1$.
In his final Example~4 in \cite{Za06} Zagier outlines the asymptotics
\begin{align}
-\ln\prod_{m=0}^\infty(1-e^{-(m+a)x})
&\sim\frac{\zeta(2)}{x}-\ln x^{\zeta(0,a)}+\ln\zeta'(0,a)
-\sum_{n=1}^\infty\frac{B_n}{n\cdot n!}\,\frac{B_{n+1}(a)}{n+1}\,(-x)^n
\nonumber\\
&=\frac{\pi^2}{6x}+\ln x^{a-\frac12}+\ln\frac{\Gamma(a)}{\sqrt{2\pi}}
-\sum_{n=1}^\infty\frac{(-1)^nB_n}{n\cdot(n+1)!}\,B_{n+1}(a)x^n
\label{eq:asym}
\end{align}
as $x\to0$, in which $\zeta(s,a)$ denotes the Hurwitz zeta function, $B_n(t)$ and $B_n=B_n(1)$ are the Bernoulli polynomials and numbers respectively.
Note that $B_n=B_n(1)=0$ for $n>1$ odd.
More generally, $B_n(1-a)=-B_n(a)$ for $n>1$ odd implying that the asymptotics
\begin{align*}
-\ln\prod_{m=0}^\infty(1-e^{-(m+a)x})(1-e^{-(m+1-a)x})
\sim\frac{\pi^2}{3x}+\ln\frac{\sqrt\pi}{\sqrt{2}\,\sin\pi a}
-\frac12\Big(a^2-a+\frac16\Big)x+o(x^N)
\end{align*}
as $x\to0$, where $N>1$ can be taken arbitrary.
In the case of rational $a$, writing it as $a/b$ with $0<a<b$ integers to the lowest terms and taking $2\pi bx$ for $x$, the asymptotics reads
\begin{align*}
-\ln(q^a,q^{b-a};q^b)_\infty\big|_{q=e^{-2\pi x}}
\sim\frac{\pi}{6bx}+\ln\frac{\sqrt\pi}{\sqrt{2}\,\sin(\pi a/b)}
-bB_2(a/b)\pi x+o(x^N)
\end{align*}
and further
\begin{align*}
q^{bB_2(a/b)/2}(q^a,q^{b-a};q^b)_\infty\big|_{q=e^{-2\pi x}}
\sim\frac{\sqrt{2}\,\sin(\pi a/b)}{\sqrt\pi}\,e^{-\pi/(6bx)}
\,\big(1+o(x^N)\big)
\end{align*}
as $x\to0$ for any $N>1$.
This agrees with the modular behaviour (hinted in Section~\ref{sec:1}) of the symmetric products $q^{bB_2(a/b)/2}(q^a,q^{b-a};q^b)_\infty$.
At the same time the asymptotics of an individual product $(q^a;q^b)_\infty$ for $a\ne b/2$ and $a\ne b$ as $q=e^{-2\pi x}\to1$ is clearly very different, since the terms for even $n$ in \eqref{eq:asym} contribute to it.
This also explains why the functions $K_4(q)$ and $K_5(q)$, even when multiplied by a rational power of~$q$, are \emph{not} modular (but, perhaps, are components of a vector-valued modular function with a controllable mock-theta-like behaviour \cite{On09,Za09} at roots of unity).
What we could check numerically was that the functions $K_4(\pm q^k)$ and $K_5(\pm q^k)$ for various $k$ (and a particular sign choice for each~$k$) do not seem to be linearly related to each other with coefficients from the field generated by $q^{k/24}$ and modular functions.

At the same time the modular function $q^{1/12}K(q)$, with $K(q)=1/(q,q^2;q^3)_\infty$, alone satisfies very simple modular equations:
\[
K(q)=\frac{(q^6;q^6)_\infty(q^8;q^8)_\infty^2}
{(q^2;q^2)_\infty^2(q^{24};q^{24})_\infty}\,\big(K(-q^4)+qK(q^{16})\big)
\]
and
\[
K(q)=\frac{(q^2;q^2)_\infty(q^3;q^3)_\infty^2(q^8;q^8)_\infty^2(q^{12};q^{12})_\infty}
{(q;q)_\infty^2(q^4;q^4)_\infty(q^6;q^6)_\infty^2(q^{24};q^{24})_\infty}\,\big(K(-q^4)-qK(q^{16})\big);
\]
these resemble the tangled equations \eqref{RR} for the Rogers--Ramanujan functions.

It seems to be appropriate to comment on how one experimentally discovers the product sides of identities like \eqref{prod}.
The principal point is that any (formal) power series $F(q)\in 1+q\mathbb Z[[q]]$ can be cast as the product $\prod_{k=1}^\infty(1-q^k)^{r_k}$ with $r_k\in\mathbb Z$. The factorisation is unique as follows from
\[
\sum_{k=1}^\infty\frac{kr_k\,q^k}{1-q^k}
=-q\,\frac{F'(q)}{F(q)};
\]
writing the latter series as $\sum_{n=1}^\infty c_nq^n$ we deduce that $c_n=\sum_{k\mid n}kr_k$, hence with the help of the M\"obius inversion formula we find out that
\[
r_k=\frac1k\sum_{d\mid k}\mu\Big(\frac kd\Big)c_d
\]
(notice that this formula obscures the integrality of the exponents~$r_k$).
This strategy makes it clear how the products like $K_1(q),\dots,K_5(q)$ for the sum sides of the Kanade--Russell conjectures can be found.
The sums in \eqref{prod} live in $1+q\mathbb Z[\omega][[q]]$; whenever we suspect a product of the type $K_6(q)$ or $K_7(q)$ for such a power series $F(q)$ we can look for one for $F(q)\ol F(q)\in1+q\mathbb Z[[q]]$ (and, in more general cases, when $F(q)\in 1+q\mathbb Z[\zeta][[q]]$ with $\zeta$ a root of unity, for the product of $F(q)$ and all its Galois conjugates) and then reconstruct the corresponding exponents for individual factors $1-\omega q^k$ and $1-\ol\omega q^k$ using the fact that the products for $F(q)$ and $\ol F(q)$ are conjugate to each other.

Since
\[
K_6(q)\ol K_6(q)=\frac{(q^6;q^9)_\infty}{(q;q^3)_\infty(q^2;q^3)_\infty^2}
\quad\text{and}\quad
K_7(q)\ol K_7(q)=\frac{(q^3;q^9)_\infty}{(q;q^3)_\infty^2(q^2;q^3)_\infty}
\]
are asymmetric, neither $K_6(q)$ nor $K_7(q)$ is modular. Furthermore,
\[
K_6(q)\ol K_7(q)=\frac{(q^3,q^6;q^9)_\infty(\omega q,\ol\omega q^2,\omega q^3,\ol\omega q^3;q^3)_\infty}{(q,q^2;q^3)_\infty}
=\frac{(\omega q;\omega q)_\infty}{(q;q)_\infty}
\]
and its conjugate $\ol K_6(q)K_7(q)$ are again not modular functions.

We conclude this section echoing \cite{Za22}: the asymptotics at roots of unity for a particular $q$-sum-to-product identity can serve, at least in principle, as a ground for its proof.
Though it looks like a doable task for the \emph{logarithm} of the product side, there seem to be no efficient strategies to make it work for the sum side.

\section{Finite identities}
\label{sec:5}

In fact, it is quite suggestive that manipulations of the type $q\mapsto\pm q^k$ are natural at the level of \emph{finite} ($q$-polynomial) versions of the sum sides of Rogers--Ramanujan(-type) identities.
These usually originate from combinatorial interpretations, and many are recorded in the literature; we limit ourselves to citing \cite{Si03} and \cite{UZ21} for a historical overview and references provided there.
In spite of this personal belief there seems to be no evidence for existence of finite versions of equations \eqref{RR} and \eqref{septic}.

\medskip\noindent
\textbf{Acknowledgements.}
I am heartily thankful to Hjalmar Rosengren, Ali Uncu and Ole Warnaar for discussions on the matter.
Part of the work was done during my visit in the Max-Planck Institute of Mathematics (Bonn) in summer 2023;
I thank the institute for excellent working conditions provided.


\begin{thebibliography}{99}

\bibitem{GIS99}
\textsc{K. Garrett}, \textsc{M.\,E.\,H. Ismail} and \textsc{D. Stanton},
Variants of the Rogers--Ramanujan identities,
\emph{Adv. Appl. Math.} \textbf{23} (1999), no.~3, 274--299.

\bibitem{Hi21}
\textsc{D. Hickerson},
\emph{Unpublished notes} (2021).

\bibitem{IS06}
\textsc{M.\,E.\,H. Ismail} and \textsc{D. Stanton},
Ramanujan continued fractions via orthogonal polynomials,
\emph{Adv. Math.} \textbf{203} (2006), no.~1, 170--193.

\bibitem{IZ18}
\textsc{M.\,E.\,H. Ismail} and \textsc{D. Stanton},
$q$-Bessel functions and Rogers--Ramanujan type identities,
\emph{Proc. Amer. Math. Soc.} \textbf{146} (2018), no.~9, 3633--3646.

\bibitem{KR15}
\textsc{S. Kanade} and \textsc{M.\,C. Russell},
\texttt{IdentityFinder} and some new identities of Rogers--Ramanujan type,
\emph{Exp. Math.} \textbf{24} (2015), no.~4, 419--423.

\bibitem{On09}
\textsc{K.~Ono},
Unearthing the visions of a master: harmonic Maass forms and number theory,
in: \emph{Current developments in mathematics} (Proceedings of the Harvard--MIT conference, 2008), pp.~347--454
(International Press, Somerville, MA, 2009).

\bibitem{Ro94}
\textsc{L.\,J. Rogers},
Second memoir on the expansion of certain infinite products,
\emph{Proc. London Math. Soc.} \textbf{25} (1894), 318--343. 

\bibitem{Ro24}
\textsc{H. Rosengren},
A new (but very nearly old) proof of the Rogers--Ramanujan identities,
\emph{SIGMA} \textbf{20} (2024), Paper no.~059, 10~pages.

\bibitem{Ro25}
\textsc{H. Rosengren},
New proofs of the septic Rogers--Ramanujan identities,
\emph{Int. J. Number Theory} (2025), to appear;
\emph{Preprint} \href{http://arxiv.org/abs/2302.02312}{\texttt{arXiv:2302.02312 [math.NT]}} (2023).

\bibitem{Si03}
\textsc{A.\,V. Sills},
Finite Rogers--Ramanujan type identities,
\emph{Electron. J. Combin.} \textbf{10} (2003), Art.~\#R13, 122~pages.

\bibitem{UZ21}
\textsc{A. Uncu} and \textsc{W. Zudilin},
Reflecting (on) the modulo 9 Kanade--Russell (conjectural) identities,
\emph{S\'em. Lothar. Combin.} \textbf{85} (2021), Art.~B85e, 17~pages;
\emph{Preprint} \href{http://arxiv.org/abs/2106.02959}{\texttt{arXiv:2106.02959 [math.NT]}} (2021).

\bibitem{Za06}
\textsc{D. Zagier},
\href{https://people.mpim-bonn.mpg.de/zagier/files/tex/MellinTransform/fulltext.pdf}%
{The Mellin transform and other useful analytic techniques},
\emph{6.7 Appendix} to E. Zeidler, \emph{Quantum Field Theory I: Basics in Mathematics and Physics. A Bridge Between Mathematicians and Physicists}
(Springer, Berlin--Heidelberg--New York, 2006), 305--323.

\bibitem{Za09}
\textsc{D.~Zagier},
Ramanujan's mock theta functions and their applications (after Zwegers and Ono--Bringmann),
in: \emph{S\'eminaire Bourbaki}, vol.~2007/2008,
\emph{Ast\'erisque} \textbf{326} (2009), Exp.~no.~986, 143--164.

\bibitem{Za22}
\textsc{D. Zagier},
Power partitions and a generalized eta transformation property,
\emph{Hardy--Ramanujan J.} \textbf{44} (2021), 1--18.

\end{thebibliography}
\end{document}